\documentclass[12pt,a4paper]{article}
\usepackage{lmodern}
\usepackage{longtable}
\usepackage{rotating}
\usepackage{url}
\usepackage{natbib}

\usepackage{xr}

\usepackage{comment}
\usepackage{mathtools}
\usepackage{array}
\usepackage{ amssymb }
\usepackage{mathrsfs}
\usepackage{cprotect}
\usepackage{arydshln}
\usepackage{ dsfont }
 \usepackage{amsthm}
\usepackage{hyperref}

\usepackage[normalem]{ulem}
\usepackage{bm}
\usepackage{enumitem}
\usepackage{graphicx}
\usepackage[multidot]{grffile}

\newcommand{\institute}[1]{#1}

\newtheorem{theorem}{Theorem}

\newtheorem{proposition}{Proposition}
\newtheorem{lemma}{Lemma}

\usepackage{filecontents}

\begin{document}


\def\spacingset#1{\renewcommand{\baselinestretch}%
{#1}\small\normalsize} \spacingset{1}


   \title{\bf Goodness-of-fit tests for functional form of Linear Mixed effects Models}
  \author{Rok Blagus and Jakob Peterlin\\
    Institute for Biostatistics and Medical Informatics, University of Ljubljana\\
   }

\date{}
 \maketitle

 \textit{\footnotesize{
 \institute{Rok Blagus
               Institute for Biostatistics and Medical Informatics, University of Ljubljana\\
               Vrazov trg 2, Ljubljana, Slovenia \\
               \href{mailto:rok.blagus@mf.uni-lj.si}{rok.blagus@mf.uni-lj.si}
               }
 }}

 \spacingset{1}
\begin{abstract}
Linear mixed effects models (LMMs) are a popular and powerful tool for analyzing clustered or repeated observations for numeric outcomes. LMMs consist of a fixed and a random component, specified in the model through their respective design matrices. Checking if the two design matrices are correctly specified is crucial since mis-specifying them can affect the validity and efficiency of the analysis. We show how to use random processes defined as cumulative sums of appropriately ordered model's residuals to test if the functional form of the fitted LMM is correctly specified. We show how these processes can be used to test goodness-of-fit of the functional form of the entire model, or only its fixed and/or random component. Inspecting plots of the proposed processes is shown to be highly informative about the potential mis-specification of the functional form of the model, providing clues for potential improvement of the model's fit. We show how the visual inspection can be objectified by using a novel procedure for estimating $p$-values which can be based on sign-flipping/bootstrap or simulations and show its validity by using theoretical results and a large Monte Carlo simulation study. The proposed methodology can be used with LMMs with multi-level or crossed random effects.
\end{abstract}



\spacingset{1} 

\section{Introduction}

Assume a (single-level) linear mixed effects model (LMM) \citep{Laird82},
\begin{equation}\label{mod}\bm{y}_i=\bm{X}_i\bm\beta+\bm{Z}_i\bm{b}_i+\bm{\epsilon}_i\mbox{, }i=1,\ldots,n \end{equation}
where $\bm{y}_i=(y_{i1},\ldots,y_{in_i})^T$, $y_{ij}$, $i=1,\ldots,n$, $j=1,\ldots,n_i$ is the  outcome for individual $j$ within cluster $i$; $\bm{X}_i$ is an $n_i\times p$ assumed fixed design matrix of the fixed effects $\bm\beta=(\beta_1,\ldots,\beta_p)^T$ with full rank, $\bm{Z}_i$ is an $n_i\times k$ assumed fixed design matrix of the random effects $\bm{b}_i=(b_{i1},\ldots,b_{ik})^T$,
$$\bm{b}_i\sim \bm{N}_k(\bm 0, \bm D), $$
$\bm D$ is positive-definite $k\times k$ matrix, and $\bm{\epsilon}_i=(\epsilon_{i1},\ldots,\epsilon_{in_i})^T$ is random noise,
$$\bm{\epsilon}_i\sim\bm{N}_{n_i}(\bm 0,\sigma^2 \bm{I}_{n_i}), $$
where $\bm{N}_d(\cdot)$ denotes a $d$-dimensional multivariate normal distribution and $\bm{I}_d$ is a $d\times d$ identity matrix. $\bm{b}_i$ and $\bm{\epsilon}_i$ are assumed independent. Let $\hat{\bm\beta}$ and $\hat{\bm{b}}_i$ denote the estimators of $\bm\beta$ and $\bm{b}_i$, respectively. Let $N=\sum_{i=1}^n n_i$ denote the total sample size. 
We assume throughout that $\sum_{i=1}^n(n_i-k)>0$ and at least one $\bm{Z}_i$ has full rank so that the model is identifiable for $\bm D$ and $\sigma^2$ \citep{demidenko2005mixed}.

Checking if the assumed LMM is correctly specified is important since model mis-specification affects the validity and efficiency of regression analysis. The most commonly used techniques for assessing the goodness-of-fit for LMMs are graphical tools such as residual plots \citep{BatesBook,wu2009mixed}. These procedures are highly subjective and often completely uninformative. \citet{Loy17} derived an approach based on the concept of visual $p$-values \citep{Mahbubul13} to make such plots less subjective. Having to rely on human experts observing the plots is however impractical. While there are numerous formal tests to check the distributional assumptions of model (\ref{mod}) \citep{jiang2001,Ritz04,Claeskens2009,Efendi17}, only few tests are available for checking its functional form. \citet{Tang14} derived a test statistic for the functional form of the fixed effects. 
The test 
involves a partition of the fixed effects design matrix. The performance of the test depends on the choice of partition and can be poor if the partition is not selected appropriately. \citet{Lee12} use a permutation approach for the inference regarding the inclusion or exclusion of the random effects in the LMMs. They show that with small to moderate samples this leads to correct inference, whereas using Wald, score, and likelihood ratio tests \citep{Self87,Stram95,Drikvandi12} does not. \citet{Pan05} propose to test the functional form of generalized LMMs (GLMMs) by considering cumulative sum (cusum) of ordered residuals. For LMMs this approach has no power against alternatives where the fixed effects design matrix is correctly specified but the random effects design matrix is not. To the best of our knowledge there is no test available for testing the functional form of the entire LMM, i.e. to test the correct specification of the fixed and random effects design matrices. We propose to base the tests for the null hypothesis
$$H_0^O: \mbox{functional form of   model (\ref{mod}) is correctly specified}, $$
by inspecting the following cusum process,
\begin{equation}
\label{proc}
W_N^O(t):=\frac{1}{\sqrt{n}}\sum_{i=1}^n  \sum_{j=1}^{n_i}e_{ij}^{C,S} I( \hat{y}_{ij}^I\leq t ),
\end{equation}
where  $I(\cdot)$ is the indicator function,
$$\hat{\bm{y}}_i^I:=(\hat{y}_{i1}^I,\ldots,\hat{y}_{in_i}^I)^T= \bm{X}_i\hat{\bm\beta}-\bm{Z}_i\hat{\bm{b}}_i$$
are \textit{individual} predicted values and
$$\bm{e}_i^{C,S}=(e_{i1}^{C,S},\ldots,e_{in_i}^{C,S})^T=\hat{\bm{S}}_i  \bm{e}_i^C $$
are transformed residuals, $\bm{e}_i^C$, weighted by some $n_i\times n_i$ weight matrix, $\hat{\bm{S}}_i$. We will show in the paper how, with a reasonably large $n$ or some assumptions on the fixed and random effects design matrices, \textit{individual} residuals,
$$\bm{e}_i^I:=(e_{i1}^I,\ldots,e_{in_i}^I)^T=\bm{y}_i-\bm{X}_i\hat{\bm\beta}-\bm{Z}_i\hat{\bm{b}}_i, $$
can be transformed in such a way that, when the functional form of the fitted LMM is correctly specified, the transformed residuals and predicted values will be uncorrelated, whereas this will not hold when the functional form is mis-specified. Under $H_0^O$ the process $W_N^O(t)$ is expected to fluctuate around zero and test statistic, $T$, can be any function mapping to the positive part of the real line such that large values give evidence for the lack-of-fit. It is essential to use \textit{individual} predicted values in the definition of $W_N^O(t)$, but \textit{individual} residuals could also be replaced with \textit{cluster} residuals,
 $$ \bm{e}_i^P:=(e_{i1}^P,\ldots,e_{in_i}^P)^T= \bm{y}_i-\bm{X}_i\hat{\bm\beta},$$
  requiring a different transformation.  

To test the null hypothesis,
$$H_0^F: \mbox{functional form of the fixed effects part of model (\ref{mod}) is correctly specified}, $$
 we define the following process,
\begin{equation}
\label{proc1}
W_N^F(t):=\frac{1}{\sqrt{n}}\sum_{i=1}^n \sum_{j=1}^{n_i}e_{ij}^{I,S} I( \hat{y}_{ij}^P\leq t ),
\end{equation}
where
$$\hat{\bm{y}}_i^P:=(\hat{y}_{i1}^P,\ldots,\hat{y}_{in_i}^P)^T=\bm{X}_i\hat{\bm\beta}$$
are \textit{cluster} predicted values and
$$\bm{e}_i^{I,S}=(e_{i1}^{I,S},\ldots,e_{in_i}^{I,S})^T=\hat{\bm{S}}_i \bm{e}_i^I,$$
are \textit{individual} residuals weighted by some $n_i\times n_i$ weight matrix $\hat{\bm{S}}_i$. Note the use of different ordering, which enables investigating only the fixed effects part of the model. We will show that when the fixed effects design matrix is correctly specified, the process $W_N^F(t)$,  is expected to fluctuate around zero even when the random effects design matrix is mis-specified. 
A similar process as (\ref{proc1}) was considered by \citet{Pan05} where (standardized) \textit{cluster} residuals are used instead of  (standardized) \textit{individual} residuals. 
We will explain how the process defined in Eq. (\ref{proc1}) can be modified to check the goodness-of-fit of any subset of the fixed effects design matrix. This was considered also by \citet{Pan05}, but when considering subsets of more than a single fixed effects covariate at a time, their approach is computationally more demanding since it requires multivariate ordering hindering also the visual presentation. 

Using random processes constructed as cumulative sum(s) of the model's residuals for goodness-of-fit testing is common for linear models (LM) \citep{Christensen15,Lin02,Stute98a,Diebolt01,Su91,Fan01,Stute98b,Blagus19} and was used also for marginal models (MMs) \citep{Lin02} and single-level GLMMs \citep{Pan05}. The challenging part for all applications of the cusum processes in this context is obtaining the null distribution of $T$. Given the complexity of the problem introduced by the dependence among the residuals asymptotic distribution for even the most trivial test statistics is analytically intractable. In LMs the null distribution of the test statistics is obtained by using bootstrap \citep{Stute98a}, simulations \citep{Su91,Lin02} or permutations \citep{Blagus19}. The simulation approach was used also for MMs \citep{Lin02} and single-level GLMMs \citep{Pan05}. We show how to use sign-flipping \citep{Winkler14}, wild bootstrap \citep{Stute98a} and simulations to correctly approximate the null distribution of the proposed processes. We also propose a novel, more powerful simulation approach. While we, for brevity of the presentation, only consider single-level LMMs, the proposed methodology can straightforwardly be extended to multiple and multi-level random effects.

Note that we are not interested in the distributional assumptions of model (\ref{mod}). Throughout the paper we refer to $H_0^O$ and $H_0^F$ simply as $H_0$ where the referenced null hypothesis should be understood from the context. The rest of the paper is organized as follows. First, we introduce some additional notation. Then we present the proposed methodology for testing the functional form of LMMs, show its asymptotic validity under the null and alternative hypotheses and showcase the finite sample performance by a selection of Monte-Carlo simulation results. An application to real data example is then given and the paper concludes with a summary of the most significant findings and possibilities for future research.



\section{Notation, definitions and estimation for (single-level) LMMs}
\label{est.lmm}

Let $\theta(\bm{X}_i)$ be the marginal mean of the correctly specified LMM,
$$ \theta(\bm{X}_i)=E(\bm{y}_i),$$
and let $\bm{\psi}_i$ be the respective marginal variance,
$$\bm{\psi}_i=var(\bm{y}_i).$$
The correctly specified model is then
$$\bm{y}_i=\theta(\bm{X}_i)+\bm{\xi}_i\mbox{, }i=1,\ldots,n, $$
where $\bm{\xi}_1,\ldots,\bm{\xi}_n$ are independent with $E(\bm{\xi}_i)=0$, $i=1,\ldots,n$, and $var(\bm{\xi}_i)=\bm{\psi}_i$.

Rewrite the assumed LMM given in (\ref{mod}) as
$$\bm{y}=\bm{X}\bm\beta+\bm{Z}\bm{b}+ \bm{\epsilon},$$
where $\bm{y}=(\bm{y}_1,\ldots,\bm{y}_n)^T$, $\bm{X}=( \bm{X}_1,\ldots,\bm{X}_n )^T$ is a $N\times p$ design matrix of the fixed effects, $\bm{Z}$, a $N\times nk$ block diagonal with $i$-th diagonal block equal to $\bm{Z}_i$, is the design matrix of the random effects, $\bm{b}=( \bm{b}_1,\ldots, \bm{b}_n)^T$ is a $nk$-vector of random effects and $\bm{\epsilon}=(\bm{\epsilon}_1,\ldots,\bm{\epsilon}_n)^T$ is a $N$-vector of random errors.

For known marginal covariance matrix
$$\bm{V}_i:=\bm{Z}_i \bm D \bm{Z}_i^T+ \sigma^2 \bm{I}_{n_i}\mbox{, }i=1,\ldots,n, $$
the estimates of the fixed effects, $\hat{\bm\beta}$, are obtained by using generalized least squares
$$\hat{\bm\beta}=  \bm{H}^{-1} \sum_{i=1}^n  \bm{X}_i^T \bm{V}_i^{-1} \bm{y}_i,$$
where
$$\bm{H}:=\sum_{i=1}^n  \bm{X}_i^T \bm{V}_i^{-1}\bm{X}_i. $$

The unique parameters in $\bm{V}_i$ are estimated either by the method of maximum likelihood (ML) or restricted ML (REML).
The estimates are then obtained by an iterative procedure, where at each step  the marginal covariance matrix is estimated   and corresponding fixed effects estimates are calculated until convergence. Based on the fitted LMM, the Best Linear Unbiased Predictors (BLUPS), $\hat{\bm{b}}_i$ are obtained from
$$\hat{\bm{b}}_i=\bm D  \bm{Z}_i^T \bm{V}_i^{-1}(\bm{y}_i-\bm{X}_i\hat{\bm\beta})=\bm D  \bm{Z}_i^T \bm{V}_i^{-1}\bm{e}_i^P, $$
where the unknown quantities are replaced by their respective estimates.

Define $\hat{\bm{G}}_i$, a consistent estimator of
\begin{equation}\label{G} \bm{G}_i:=\bm{I}_{n_i}-\bm{Z}_i \bm D  \bm{Z}_i^T \bm{V}_i^{-1}=\sigma^2\bm{V}_i^{-1},\end{equation}
which is obtained by replacing the unknown quantities in (\ref{G}) with their respective consistent estimators, $\hat{\bm D}$ and  $\hat{\sigma}^2$. Observe that $\hat{\bm{G}}_i$ is symmetric and that, when the random effects are estimated by using the (estimated) BLUPS,
\begin{equation}
\label{ind.to.clust}
\bm{e}_i^I=\hat{\bm{G}}_i\bm{e}_i^P=\sigma^2\hat{\bm{V}}_i^{-1}\bm{e}_i^P.
\end{equation}


Transformed residuals used in equation (\ref{proc}) are defined as
\begin{equation}
\label{eq:rot}
\bm{e}^C:=(\bm{e}_1^C,\ldots,\bm{e}_n^C)^T=\bm{e}^I - \hat{\bm{A}}\hat{\bm{B}}^{+}    \bm{Z}\hat{\bm{b}},
\end{equation}
where $\bm{e}^I=( \bm{e}_1^I,\ldots,\bm{e}_n^I )^T$, $\hat{\bm{b}}=( \hat{\bm{b}}_1,\ldots, \hat{\bm{b}}_n)^T$ and $\hat{\bm{A}}$ and $\hat{\bm{B}}$ are respective consistent estimators of block matrices
\begin{small}
\begin{equation}
\label{eq:A}
\bm{A}=cov(\bm{e}^I,\hat{\bm{y}}^I )=\sigma^2 \bm{V}^{-1}\left[ \bm{V}- \bm{X}  \bm{H}^{-1}   \bm{X}^T \right]\left( \bm{I}_N-\sigma^2 \bm{V}^{-1}  \right),
\end{equation}
\end{small}
and
\begin{small}
\begin{equation}
\label{eq:B}
\bm{B}=var(\hat{\bm{y}}^I)=\left( \bm{I}_N-\sigma^2 \bm{V}^{-1}  \right)\left[ \bm{V}- \bm{X}  \bm{H}^{-1}   \bm{X}^T \right]\left( \bm{I}_N-\sigma^2 \bm{V}^{-1}  \right),
\end{equation}
\end{small}
where $\hat{\bm{y}}^I=(\hat{\bm{y}}^I_1,\ldots,\hat{\bm{y}}^I_n)^T$, $\bm{V}$ is $N\times N$ block diagonal matrix with $i$-th diagonal block equal to $\bm{V}_i$ and where $\hat{\bm{B}}^{+}$  denotes the Moore-Penrose inverse of $\hat{\bm{B}}$ (observe that $\hat{\bm{B}}$ is not invertible).

Define
\begin{equation}
\label{eq:Jot}
\bm{J}:=\sigma^2 \bm{V}^{-1}- \bm{A}\bm{B}^+\bm{Z}\left( \bm{I}_n\otimes \bm{D}  \right)\bm{Z}^T \bm{V}^{-1},
\end{equation}
where $\otimes$ is the Kronecker product. Then, using (\ref{ind.to.clust}), we obtain
$$\bm{e}^C= \hat{\bm{J}}\bm{e}^P,$$
where $\bm{e}^P=(\bm{e}^P_1,\ldots,\bm{e}^P_n)^T$ and $\hat{\bm{J}}$ is a consistent estimator of $\bm{J}$ obtained replacing the unknown quantities in (\ref{eq:Jot}) with their respective consistent estimators.

Define
\begin{equation}
\label{eq:chi}
\bm{\chi}_i(z,\bm\beta,\bm{b}_i):=(I( \bm{X}_{i1}\bm\beta+\bm{Z}_{i1}\bm{b}_i \leq z),\ldots,I( \bm{X}_{in_i}\bm\beta+\bm{Z}_{in_i}\bm{b}_i \leq z))^T=\bm{I}(\bm{X}_{i}\bm\beta+\bm{Z}_{i}\bm{b}_i \leq z\bm{1}),
\end{equation}
where $\bm{X}_{ij}$ and $\bm{Z}_{ij}$, $i=1,\ldots,n$, $j=1,\ldots,n_i$ are the $j$-th rows of the fixed and random effects design matrices, respectively and $\bm{1}=(1,\ldots,1)^T$ is a $n_i$-vector of ones. Let $\bm{\chi}(z,\bm\beta,\bm{b})=(\bm{\chi}_1(z,\bm\beta,\bm{b}_1),\ldots,\bm{\chi}_n(z,\bm\beta,\bm{b}_n))^T$.



Recall that the residuals are standardized by some weight matrix, $\hat{\bm{S}}_i$, which if properly defined can increase the convergence rate and power \citep{Pan05}. Let $\bm{S}$ be a $N\times N$ block diagonal matrix, with $i$-th diagonal block equal to $\bm{S}_i$ and let $\hat{\bm{S}}$ be a consistent estimator of $\bm{S}$. We used $\hat{\bm{S}}_i=\hat{\bm{V}}_i^{-1/2}$; other definitions  are also possible but were not investigate here in more detail (see supplementary information).



\section{Various aspects of goodness-of-fit of the functional form for fitted LMMs}
\label{gof.lmm}

There are three  aspects of goodness-of-fit of the functional form for LMMs. The first is the correct specification of the functional form of the entire fitted model. Are the assumed design matrices for the fixed and random effects correctly specified? This can be answered by inspecting  $W_N^O(t)$.  Second, given the assumed random effects structure, is the assumed  design matrix for the fixed effects correctly specified? This can be answered by inspecting $W_N^F(t)$. And third, given the assumed design matrix for the fixed effects, is the assumed  design matrix for the random effects correctly specified? Importantly, it turns out that it is not necessary to answer the last question directly, which in fact would be very difficult. We formally prove that when the fixed effects design matrix is correctly specified, the test based on the process $W_N^F(t)$   will be robust against mis-specification of the random effect design matrix. Intuitively, this should hold since the estimator of $\bm\beta$ is consistent given only the correct specification of the fixed effect design matrix \citep{Zeger88}. Therefore, when there is enough evidence to reject the goodness-of-fit of the entire model   but at the same time there is no evidence against the goodness-of-fit of the fixed effects part of the model, this implies mis-specification of the random effects part of the model.

The proposed processes, $W_N^O(t)$ and $W_N^F(t)$, can be visualized by plotting them against $t$ from were it is possible to (subjectively) assess the adequacy of the fitted model by comparing the observed process with random realizations of the processes under $H_0$ and can be objectified by calculating  $p$-values for some reasonable test statistics. We show later how this can be obtained by using sign-flipping/bootstrap or simulations.

For the second question the process can be modified by defining
$$W_N^{F^S}(t):=\frac{1}{\sqrt{n}}\sum_{i=1}^n  \sum_{j=1}^{n_i}e_{ij}^{I,S} I( \sum_{l} X_{ij,l}\hat\beta_l  \leq t ), $$
where $X_{ij,l}$ is the $l$-th column, $l\subset\{1,\ldots,p\}$, and $j$-th row, $j=1,\ldots,n_i$, of the design matrix $\bm{X}_i$, $i=1,\ldots,n$, and the sum extends only over some subset of the column's  of the fixed effects design matrix. When using $W^{F^S}(t)$ one tests for a possible lack-of-fit which is only due to the specified fixed effects covariates (by inspecting $W^{F^S}(t)$ it is also possible to detect an omission of an important interaction effect, see \citet{Blagus19}). All the processes can be modified by taking the sum only within a window specified by some positive constant $c>0$ as proposed by \citet{Lin02}. E.g. in (\ref{proc}) one could use $I(t-c<\hat{y}_{ij}^I\leq t)$. 
The constructed random processes tend to be dominated by the residuals with small fitted  values, which can potentially lead to a loss of power \citep{Lin02}. While specifying a meaningful value of $c$ could potentially improve the power of the proposed tests, this was not considered here. 

The following (iterative) approach is then proposed to check if the functional form of the fitted LMM is correctly specified. First, plot $W_N^O(t)$  and $W_N^F(t)$ (and if necessary $W^{F^S}(t)$) along with their null realizations and obtain $p$-values. If none of the $p$-values is significant and the observed processes are similar as their null realizations, there is not enough evidence against the goodness-of-fit. If only the $p$-value based on $W_N^O(t)$ is significant, this indicates the mis-specification of the assumed form of the random effects. In any other case, fixed and/or  random effects design matrices might be mis-specified. In the last case, first correct the fixed effects part of the model, i.e. find the model for which the $p$-value based on  $W^F(t)$ is no longer significant. Here one can use the process $W^{F^S}(t)$ to detect the covariate (or a set of them) which is causing the lack-of-fit. The result of \citet{Lin02} showing the representative cusum processes when mis-specifying the functional form of the fixed effects in specific way, can be helpful and the $p$-values can be complemented also by using directional tests, e.g. Wald test. When there is no more evidence against the lack-of-fit for the assumed fixed effects part of the model, check $W_N^O(t)$. If it still shows a significant decline from the assumed model, this implies mis-specification of the random effects design matrix. Here the $p$-values based on the process $W_N^O(t)$ can be complemented by using directional tests, e.g. the tests proposed by \citet{Lee12}, to help specifying the form of the random effects design matrix until there is no more evidence for the lack-of-fit of the entire model.



\begin{figure}[h!]
\caption{Different processes for the simulated example; panel (A): $W_N^O(t)$ for model (\ref{mod.ex.1}), panel (B): $W_N^F(t)$ for model (\ref{mod.ex.1}), panel (C): $W_N^O(t)$ for model (\ref{mod.ex.2}), panel (D): $W_N^F(t)$ for model (\ref{mod.ex.2}), panel (E): $W_N^O(t)$ for model (\ref{mod.ex.3}) and panel (F): $W_N^F(t)$ for model (\ref{mod.ex.3}). Gray lines are the processes obtained by using 500 random sign-flips. The $p$-values are for the KS and CvM test (the smallest possible estimated $p$-value was $1/501=0.002$).   }
\label{figIll}
\begin{center}
\resizebox{120mm}{!}{\includegraphics{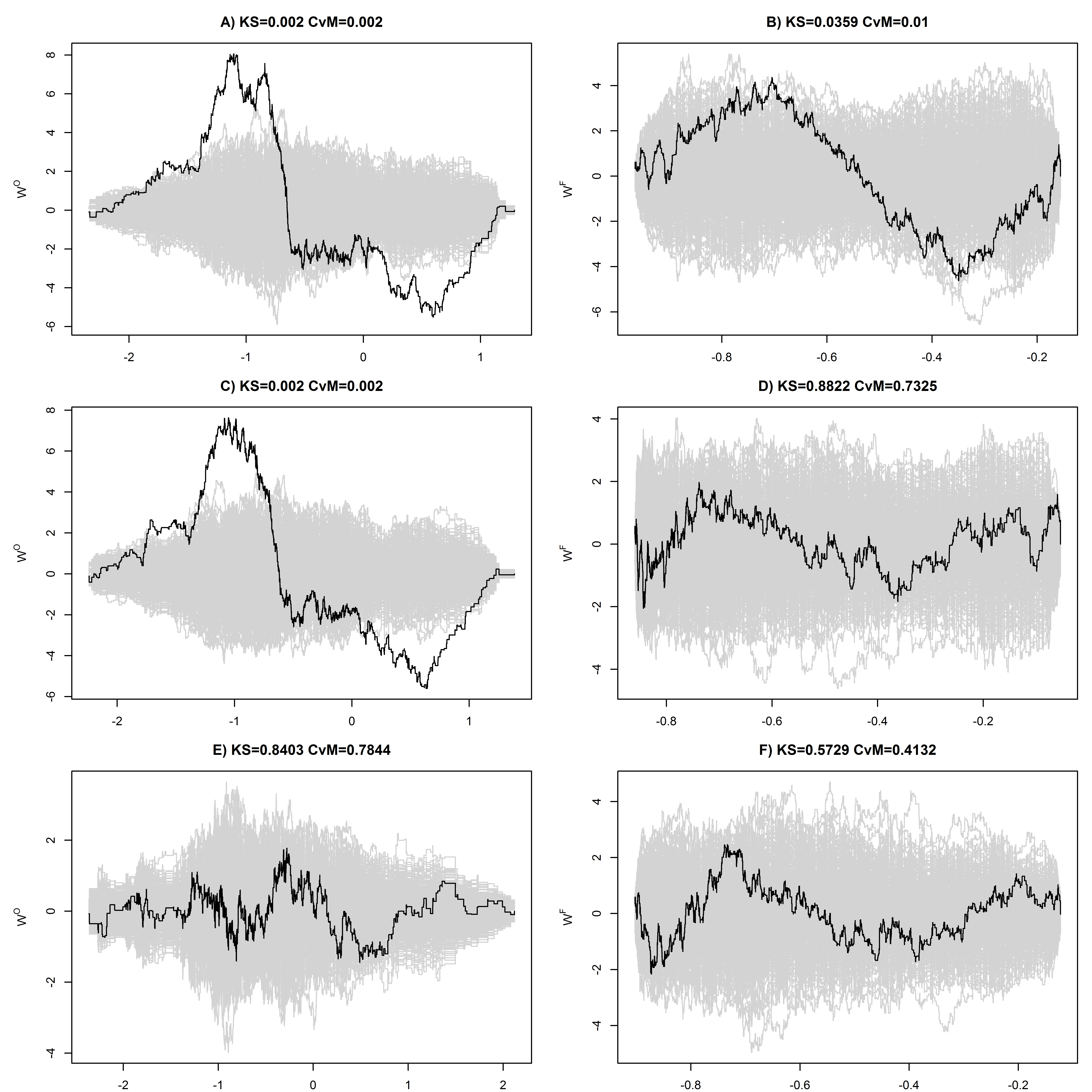}}\end{center}
\end{figure}

We illustrate the proposed approach on a simulated example where the outcome was simulated from
\begin{equation}\label{mod.ex.0} y_{ij}= - 1 +0.25 X_{ij}+0.5X_{ij}^2+b_{i,0}+b_{i,1}X_{ij}+\epsilon_{ij}\mbox{, }j=1,\ldots,n_i\mbox{, }i=1,\ldots,n,\end{equation}
where $X_{ij}\sim U(0,1)$, $\epsilon_{ij}\sim N(0,0.25)$, $b_{i,0}\sim N(0,0.25)$ and $b_{i,1}\sim N(0,1)$. $b_{i,0}$ and $b_{i,1}$ were simulated independently and $\bm{\epsilon}_i=(\epsilon_{i1},\ldots,\epsilon_{in_i})^T$ and $\bm{b}_i=(b_{i,0},b_{i,1})^T$ were also independent. There were $n=50$ clusters and for each there were $n_i=25$ subjects. Different processes are shown in Figure \ref{figIll} (black lines); gray lines are 500 random processes, obtained when using the proposed procedure based on sign-flipping (it is explained in section \ref{gof} how to obtain the gray curves and the $p$-values). 

First, we fit the model with the fixed effect covariate and a random intercept,
\begin{equation}\label{mod.ex.1} y_{ij}=  \hat\beta_0 +\hat\beta_1 X_{ij}+\hat{b}_{i,0}.\end{equation}
In Figure \ref{figIll}, panel (B) we show the process $W_N^F(t)$ for model (\ref{mod.ex.1}) from where we can with the help of results presented in \citet{Lin02},  deduce that a quadratic effect of the fixed effect covariate should be included in the model. Hence we next fit the following model,
\begin{equation}\label{mod.ex.2} y_{ij}=  \hat\beta_0 +\hat\beta_1 X_{ij}+\hat\beta_2 X_{ij}^2+\hat{b}_{i,0}.\end{equation}
The process $W_N^F(t)$ for model (\ref{mod.ex.2}) is shown as a black curve in Figure \ref{figIll}, panel (D). We can conclude that the fit of the fixed effects part of the model is satisfactory. The process $W_N^O(t)$ for model (\ref{mod.ex.2}) (Figure \ref{figIll}, panel (C)) still shows a significant decline from the assumed model, and we can deduce that the random effect part of the model is not correctly specified. We then fit the model where we add a random slope,
\begin{equation}\label{mod.ex.3} y_{ij}=  \hat\beta_0 +\hat\beta_1 X_{ij}+\hat\beta_2 X_{ij}^2+\hat{b}_{i,0}+\hat{b}_{i,1}X_{ij}.\end{equation}
The processes $W_N^O(t)$ and $W_N^F(t)$ for model (\ref{mod.ex.3}), shown in figure \ref{figIll}, panels (E) and (F), respectively, show no decline from the assumed model, hence we deduce that the fit of the model (\ref{mod.ex.3}) is satisfactory. 

When considering only $W_N^O(t)$ it is difficult to distinguish whether the lack-of-fit was due to fixed or random effects part of the model (or both). E.g. in the example the plots in Figure \ref{figIll}, panels (A) and (C) are very similar. 
Hence it is crucial to consider $W_N^O(t)$ in combination with $W_N^F(t)$ and make conclusions based on the former when the later shows no deviations for the assumed fixed effects part of the model. 
In supplementary information we show the standard plot of (standardized) residuals versus fitted values for models (\ref{mod.ex.1}), (\ref{mod.ex.2}) and (\ref{mod.ex.3}), a default \textbf{plot} option in the R's \textbf{nlme} package \citep{BatesBook}. It would take an experienced eye, or some imagination, to be able to distinguish between different plots and to decide which plot actually belongs to the correctly specified model.

\section{Approximating the null distributions of the proposed processes  }
\label{gof}

Here we show how to approximate the null distributions of the proposed processes and how to estimate $p$-values. Define
$$W_N^{O,m}(t):=\frac{1}{\sqrt{n}}\sum_{i=1}^n   \sum_{j=1}^{n_i}e_{ij}^{C,S,m} I( \hat{y}_{ij}^{I}\leq t ) \mbox{, }m=1,\ldots,M $$
$$W_N^{F,m}(t):=\frac{1}{\sqrt{n}}\sum_{i=1}^n \sum_{j=1}^{n_i}e_{ij}^{I,S,m} I( \hat{y}_{ij}^{P}\leq t )\mbox{, }m=1,\ldots,M $$
$$W_N^{F^S,m}(t):=\frac{1}{\sqrt{n}}\sum_{i=1}^n  \sum_{j=1}^{n_i}e_{ij}^{I,S,m} I( \sum_{l} X_{ij,l}\hat\beta_l\leq t )\mbox{, }m=1,\ldots,M $$
where $e_{ij}^{C,S,m}$ and $e_{ij}^{I,S,m}$ are obtained by fitting a LMM to $\bm{y}_i^m$, $m=1,\ldots,M$, where
$$\bm{y}_i^m=\hat{\bm{y}}_i^P+\hat{\bm{L}}_i\bm{\Pi}_i\hat{\bm{L}}_i^{-1}\bm{e}_i^P, $$
where $\hat{\bm{L}}_i$, such that $\hat{\bm{V}}_i=\hat{\bm{L}}_i\hat{\bm{L}}_i^T$, is found using Cholesky decomposition and where $\bm\Pi_i$ can be any $n_i\times n_i$ random matrix, such that $E(\bm\Pi_i)=\bm{0}$ and $var(\bm\Pi_i)=\bm{I}_{n_i}$. E.g., $\bm\Pi_i$ can be the sign flipping matrix \citep{Winkler14} as used in this paper, but it could also be a matrix constructed for wild bootstrap procedure as considered by \citet{Stute98a} for LMs or a diagonal matrix with random standard normal deviates on the diagonal. 

We estimate $p$-values by
$$\frac{1}{M+1}\sum_{m=1}^M \left(I( T^m\geq   T)+1\right), $$
where $T=g( W_N^{\bullet}(t) )$  and $T^m=g( W_N^{\bullet,m}(t) )$ for the original processes and their null approximations, respectively and where $g(\cdot)$ is some function mapping to the positive part of the real line so that large values give evidence against $H_0$.


The $p$-values for the Kolmogorov-Smirnov (KS) and Cramer-von Mises (CvM) type test statistics for the overall test which were considered here, were defined as
$$T^K:= \max_{t}|W_N^{\bullet}(t)|\mbox{, } T^{K,m}:=  \max_{t}|W_N^{\bullet,m}(t)|    \mbox{, }m=1,\ldots,M $$
and
$$T^C:=  \sum_{t}\left(W_N^{\bullet}(t)\right)^2\mbox{, }T^{C,m}:= \sum_{t}\left(W_N^{\bullet,m}(t)\right)^2\mbox{, }m=1,\ldots,M.$$

The use of the proposed test statistics is justified for $W_N^F(t)$ ($W_N^{F^S}(t)$) after noting that under $H_0$ the residuals and \textit{cluster} predicted values are uncorrelated,
$$cov(\bm{e}_i^I,\hat{\bm{y}}_j^P)=cov(\bm{e}_i^P,\hat{\bm{y}}_j^P)=\bm{0}\mbox{, }i=1,\ldots,n\mbox{, }j=1,\ldots,n ,$$
thence the process is expected to fluctuate around zero when $H_0$ holds.
To show that this also holds for $W_N^O(t)$ is more arduous. In general the residuals and \textit{individual} predicted values are correlated
$$cov(\bm{e}_i^I,\hat{\bm{y}}_i^I)\neq \bm{0}\mbox{ and  }cov(\bm{e}_i^P,\hat{\bm{y}}_i^I)\neq \bm{0}\mbox{, }i=1,\ldots,n , $$
thence using untransformed \textit{individual} (or \textit{cluster}) residuals  when constructing the process $W_N^O(t)$ would yield a process which even under $H_0$ does not fluctuate around zero. Therefore, we transform the residuals (see Eq. (\ref{eq:rot})) in such a way that under $H_0$, they are uncorrelated with the \textit{individual} predicted values (we show in the supplementary information an alternative approach where we do not transforms the residuals and construct the process directly from \textit{individual} or \textit{cluster} residuals but we show that this approach lacks a nice visual presentation and can be conservative).  To show that this holds, we first prove a general result for the Moore-Penrose pseudo inverse which can be applied to symmetric matrices, and then specialize it to our example (the proofs are given in the supplementary information).

\begin{lemma}
\label{jakob.lema.the}
Let $\bm{P},\bm{M}\in \mathds{R}^{q\times q}$ be some symmetric matrices. Let $\bm{C}=\bm{M}\bm{P}^T$ and $\bm{Q}=\bm{P}\bm{M}\bm{P}^T$ and let $\bm{Q}^+$ be the Moore-Penrose pseudo inverse of matrix $\bm{Q}$. Let $\mbox{Im}$ denote the image of a matrix transformation.  Assume that for each vector $\bm{v}\in\mbox{Im}(\bm{P})$,
$$\bm{M}\bm{v}\in \mbox{Im}(\bm{P}), $$
holds. Then
$$ \bm{C}-\bm{C}\bm{Q}^+\bm{Q}=\bm{0},$$
where $\bm{0}$ is a $q\times q$ matrix of zeros.
\end{lemma}

\begin{theorem}
\label{lemma.jakob.1}
Let $\bm{A}$ and $\bm{B}$ be as defined in (\ref{eq:A}) and (\ref{eq:B}), respectively. The equality
\begin{equation}\label{equality}\bm{A}-\bm{A}\bm{B}^+\bm{B}=\bm{0},\end{equation}
where $\bm{0}$ is a $N\times N$ matrix of zeros, holds if any of the following conditions are satisfied
\begin{itemize}
\item[(i)] matrix $\bm{B}$ is invertible;
\item[(ii)] for each $i$, $\mbox{Im}(\bm{X}_i)\subset\mbox{Im}(\bm{Z}_i)$;
\item[(iii)] $\bm{Z}_i=(1,\ldots,1)^T$, $i=1,\ldots,n$ and matrix $\bm{X}$ has a column of ones and all other columns are orthogonal with respect to the column of ones.
\end{itemize}
\end{theorem}


Due to identifiability constraint, we know that in our setting, condition (i) cannot hold. Condition (ii) holds when for $i=1,\ldots,n$, the design matrix $\bm{Z}_i$ includes at least all columns of the design matrix $\bm{X}_i$. Condition (iii) can be applied to models with only a random intercept. Defining
$$ \tilde{\bm{e}}^C:=   \bm{e}^I - \bm{A}\bm{B}^{+}  \bm{Z}\hat{\bm{b}} , $$
and using that
$$cov(\bm{e}^I,\hat{\bm{y}}^I)=cov(\bm{e}^I,\bm{Z}\hat{\bm{b}})=\bm{A}\mbox{ and } var( \bm{Z}\hat{\bm{b}} )=\bm{B},$$
it is not difficult to see that assuming that any condition of Theorem \ref{lemma.jakob.1} is satisfied, under $H_0$, the transformed residuals and \textit{individual} predicted values are uncorrelated
$$ cov(\tilde{\bm{e}}^C, \hat{\bm{y} }^I)= \bm{A}- \bm{A}\bm{B}^+\bm{B}=\bm{0}. $$
The above result also holds when defining the process based on transformed \textit{cluster} residuals (see supplementary information). The next results will show, however, that when $n\rightarrow\infty$ equality (\ref{equality}) applies for general design matrices  $\bm{X}$ and $\bm{Z}$, thus showing that for a large $n$ the process $W_N^O(t)$ is, under $H_0$, expected to fluctuate around zero also for general design matrices $\bm{X}$ and $\bm{Z}$. 


\begin{lemma}
\label{Rok}
Assuming that there exists a number $ \alpha < \infty $ such that $0<||\bm{H}||_2<\alpha$, then in the limit, when $n\rightarrow \infty$, $\bm{A}$ and $\bm{B}$ are block diagonal, with the $i$-th diagonal elements being equal to
\begin{small}
$$\tilde{\bm{A}}_i= \sigma^2 \bm{V}_i^{-1}\bm{Z}_i\bm{D}\bm{Z}_i^T $$
\end{small}
and
\begin{small}
 $$\tilde{\bm{B}}_i=  \bm{Z}_i\bm{D}\bm{Z}_i^T \bm{V}_i^{-1} \bm{Z}_i\bm{D}\bm{Z}_i^T  .$$
\end{small}
\end{lemma}

\begin{theorem}
\label{jes}
Let $\tilde{\bm{A}}_i$ and $\tilde{\bm{B}}_i$ be as given in  Lemma \ref{Rok}. Then
$$\tilde{\bm{A}}_i-\tilde{\bm{A}}_i \tilde{\bm{B}}_i^+\tilde{\bm{B}}_i=\bm{0},$$
where $\bm{0}$ is a $n_i\times n_i$ matrix of zeros.
\end{theorem}


Observing that in the limit, when $n\rightarrow\infty$, $\bm{A}$ and $\bm{B}$ are block diagonal, allows an alternative, computationally more efficient definition of the process $W_N^O(t)$ which according to Theorem \ref{jes} will, for a large $n$,  under $H_0$, also fluctuate around zero,
\begin{equation}
\label{proc11}
\tilde{W}_N^O(t):=\frac{1}{\sqrt{n}}\sum_{i=1}^n  \sum_{j=1}^{n_i}e_{ij}^{\tilde{C},S} I( \hat{y}_{ij}^I\leq t ),
\end{equation}
$$\bm{e}_i^{\tilde{C},S}=(e_{i1}^{\tilde{C},S},\ldots,e_{in_i}^{\tilde{C},S})^T=\hat{\bm{S}}_i  \bm{e}_i^{\tilde{C}}, $$
where transformed \textit{individual} residuals are now obtained from
$$ \bm{e}_i^{\tilde{C}}:=(e_{i1}^{\tilde{C}},\ldots,e_{in_i}^{\tilde{C}})^T=\bm{e}_i^I - \hat{\bm{A}}_i\hat{\bm{B}}_i^{+}    \bm{Z}_i\hat{\bm{b}}_i =\hat{\bm{J}}_i\bm{e}_i^P  ,  $$
where $\hat{\bm{A}}_i$ and $\hat{\bm{B}}_i$ are respective consistent estimators of the $i$-th diagonal blocks of $\bm{A}$ and $\bm{B}$, $\bm{A}_i$ and $\bm{B}_i$, respectively and where $\hat{\bm{J}}_i$ is a consistent estimator of the $i$-th diagonal block of $\bm{J}$ defined in (\ref{eq:Jot}).  

Processes  $\tilde{W}_N^{O,m}(t)$, $m=1,\ldots,M$ are defined similarly replacing $e_{ij}^{\tilde{C},S}$ with $e_{ij}^{\tilde{C},S,m}$ where the later are obtained by refitting the LMM to $\bm{y}_i^m$, re-estimating the matrices $\bm{A}_i$ and $\bm{B}_i$ for each $m$. Observe that the limit expressions, as $n\rightarrow\infty$, of matrices $\bm{A}_i$ and $\bm{B}_i$ are the same as given for matrices $\tilde{\bm{A}}_i$ and $\tilde{\bm{B}}_i$  as defined in Lemma \ref{Rok}, hence as $n\rightarrow\infty$, the estimates of $\tilde{\bm{A}}_i$ and $\tilde{\bm{B}}_i$ could equivalently be used instead and re-estimating them for each $m$ would not be necessary. However, our simulation results show that using the estimates of $\bm{A}_i$ and $\bm{B}_i$ and re-estimation for each $m$ leads to faster convergence rate  (data not shown). 

Refitting the LMM for each $m$ can however be computationally demanding with large $n$ and/or complex models. It is possible to define the processes so that refitting the LMM in each step $m=1,\ldots,M$ is not necessary. Let
$$\bm{e}(M_1,\ldots,M_n)=(\bm{e}_1(M_1),\ldots,\bm{e}_n(M_n) )^T, $$
where
$$\bm{e}_i(M_i)=\bm{e}^P_i M_i-\bm{X}_i \hat{\bm{H}}^{-1}\sum_{i=1}^n \bm{X}_i^T\hat{\bm{V}}_i^{-1}\bm{e}^P_i M_i, $$
where $M_1,\ldots,M_n$ are independent standard normal random variables as used by \citet{Pan05} for single-level GLMMs and considered in this paper, but they could also be independent and identically distributed such that $E(M_i)=0$, $var(M_i)=1$ and $|M_i|\leq c<\infty$, for some finite $c$, obtained for example by using sign-flipping or the wild bootstrap.

Then we define
\begin{equation*}
\label{eq:O.nrf}
W_N^{\tilde{O},m}(t):=\frac{1}{\sqrt{n}}\bm\chi(t,\hat{\bm{\beta}},\hat{\bm{b}})^T\hat{\bm{S}}\hat{\bm{J}}\bm{e}(M_1,\ldots,M_n),
\end{equation*}
\begin{equation}
\label{eq:O.tilde}
\tilde{W}_N^{\tilde{O},m}(t):=\frac{1}{\sqrt{n}}\sum_{i=1}^n \bm\chi_i(t,\hat{\bm{\beta}},\hat{\bm{b}}_i)^T\hat{\bm{S}}_i\hat{\bm{J}}_i\bm{e}_i(M_i),
\end{equation}
and
\begin{equation}
\label{eq:F.tilde}
W_N^{\tilde{F},m} (t):=\frac{1}{\sqrt{n}}\sum_{i=1}^n   \bm{\chi}_i(t,\hat{\bm\beta},\bm{0})^T\hat{\bm{S}}_i\hat{\bm{G}}_i\bm{e}_i(M_i)    ,
\end{equation}
where $\bm\chi_i(\cdot)$ is defined in Eq. (\ref{eq:chi}) and $\hat{\bm{J}}_i$ is the estimator of the $i$-th diagonal block of  $\bm{J}$ defined in (\ref{eq:Jot}) obtained from the original fit. Note that for each $i$ the same realization of the random variable $M_i$ is used for all $j=1,\ldots,n_i$.  Potentially more powerful approach is to define $$\bm{e}(\bm\Pi_1,\ldots,\bm\Pi_n)=(\bm{e}_1(\bm\Pi),\ldots,\bm{e}_n(\bm\Pi) )^T, $$
where
$$\bm{e}_i(\bm\Pi_i)=\hat{\bm{L}}_i \bm{\Pi}_i \hat{\bm{L}}_i^{-1}\bm{e}_i^P-\bm{X}_i \hat{\bm{H}}^{-1}\sum_{i=1}^n \bm{X}_i^T\hat{\bm{V}}_i^{-1}\hat{\bm{L}}_i \bm{\Pi}_i \hat{\bm{L}}_i^{-1}\bm{e}_i^P, $$
where $\bm{\Pi}_i$ is a random $n_i\times n_i$ matrix as defined previously. The processes are then defined as
\begin{equation*}
\label{eq:O.hat.nrf}
W_N^{\hat{O},m}(t):=\frac{1}{\sqrt{n}}\bm\chi(t,\hat{\bm{\beta}},\hat{\bm{b}})^T\hat{\bm{S}}\hat{\bm{J}}\bm{e}(\bm\Pi_1,\ldots,\bm\Pi_n),
\end{equation*}
\begin{equation}
\label{eq:O.hat}
\tilde{W}_N^{\hat{O},m}(t)=\frac{1}{\sqrt{n}}\sum_{i=1}^n \bm\chi_i(t,\hat{\bm{\beta}},\hat{\bm{b}}_i)^T\hat{\bm{S}}_i\hat{\bm{J}}_i\bm{e}_i(\bm\Pi_i),
\end{equation}
and
\begin{equation}
\label{eq:F.hat}
W_N^{\hat{F},m} (t)=\frac{1}{\sqrt{n}}\sum_{i=1}^n   \bm{\chi}_i(t,\hat{\bm\beta},\bm{0})^T\hat{\bm{S}}_i\hat{\bm{G}}_i\bm{e}_i(\bm\Pi_i).
\end{equation}
For the process $W_N^{F^S,m}(t)$  similar definitions as for $W_N^{F,m}(t)$ are used by appropriately specifying $\bm{\chi}_i(z,\bm\beta,\bm{b}_i )$. The $p$-values would then be estimated in exactly the same manner as described before. We refer to the approach where the LMM is refitted in each step $m=1,\ldots,M$ as the sign-flipping/bootstrap approach and as (novel) simulation approach where the LMM is not refitted using (a different) the same realization of the simulated random variable within each cluster.

We would proceed similarly for two-(or more-)level LMMs. The methodology developed here can be used also for crossed random effects structures by using the fact that they can alternatively be represented as the random-effects structure corresponding to a single-level model (see \citet{BatesBook}, page 165 for more details).  In the supplementary information we illustrate an extension to two-level LMMs.

\section{Asymptotic convergence of the cusum random processes under $H_0$}
\label{theor}


Here we show the asymptotic convergence of the proposed processes under $H_0$, when for some fixed and finite $n_i$, $i=1,\ldots,n$, $n\rightarrow\infty$. We will assume that $\bm{D}$ is a known positive-definite matrix and $\sigma^2$ a known constant. Asymptotic convergence of a similar process as $W_N^F(t)$ was studied by \citet{Pan05}, also in a more general setting where $\bm{D}$ and $\sigma^2$ are consistently estimated. By following the arguments presented in \citet{Pan05}, the results presented here could be extended to the setting with unknown $\bm{D}$ and $\sigma^2$ with some additional notational inconvenience. Note that by assuming known $\bm{D}$ and $\sigma^2$, the identifiability constraint is not required here.  We consider in more detail the convergence of $\tilde{W}_N^O(t)$  and $W_N^F(t)$  along with their respective null approximations, $\tilde{W}_N^{O,m}(t)$ , $\tilde{W}_N^{\tilde{O},m}(t)$ and $\tilde{W}_N^{\hat{O},m}(t)$ and $W_N^{F,m}(t)$, $W_N^{\tilde{F},m}(t)$ and $W_N^{\hat{F},m}(t)$.  Using the arguments presented in \citet{Pan05} and Lemma \ref{Rok}, convergence of $W_N^O(t)$ follows from the convergence of $\tilde{W}_N^O(t)$ (by a similar argument  convergence of $W_N^{O,m}(t)$, $W_N^{\tilde{O},m}(t)$ and $W_N^{\hat{O},m}(t)$ follow from the convergence of $\tilde{W}_N^{O,m}(t)$, $\tilde{W}_N^{\tilde{O},m}(t)$ and $\tilde{W}_N^{\hat{O},m}(t)$, respectively). Convergence of the process $W_N^{F^S}(t)$ can be established similarly as for $W_N^F(t)$ by appropriately specifying $\bm{\chi}_i(t,\bm\beta,\bm{b}_i )$, but is not included here for brevity. To facilitate theoretical investigation, we extend the processes over the entire real line, so that they are elements in the Skorokhod-space $D[-\infty,\infty]$. Some assumptions are made in order to assure stochastic equicontinuity of the processes. We also make some assumptions on the norms of the fixed matrices and moments of random vectors. For brevity the assumptions and proofs are stated and discussed in the supplementary information.

The following theorem establishes the convergence of $W_N^F(t)$, see \citet{Pan05} for an alternative proof.

\begin{theorem}
\label{theorem1}
Under $H_0$ and assumptions stated in the supplementary information, with probability one, the process $W_N^F(t)$ converges in distribution to a zero-mean Gaussian process $G_\infty$ in the Skorokhod-space $D[-\infty,\infty]$,  where  the covariance function of $G_\infty$ is
  $$K(t,s):=cov( G_{\infty}(t),G_{\infty}(s) )=\lim_{n\rightarrow \infty}K_N(t,s),$$
 where
\begin{small}
\begin{eqnarray*}
K_N(t,s)&:=&\frac{1}{n }\left[ \sum_{i=1}^n \bm{\chi}_i(t)^T \bm{S}_i \bm{G}_i\bm{V}_i \bm{G}_i\bm{S}_i^T \bm{\chi}_i(s)  \right.\\
&-&\left. \left( \sum_{i=1}^n  \bm{\chi}_i(t)^T  \bm{S}_i \bm{G}_i\bm{X}_i \right)  \left( \sum_{i=1}^n \bm{X}_i^T  \bm{V}_i^{-1}   \bm{X}_i\right)^{-1}\left( \sum_{i=1}^n  \bm{X}_i^T \bm{G}_i\bm{S}_i^T\bm{\chi}_i(s) \right)  \right],\\
\end{eqnarray*}
\end{small}
where
$$\bm{\chi}_i(z):=(I( \bm{X}_{i1}\bm\beta \leq z),\ldots,I( \bm{X}_{in_i}\bm\beta \leq z))^T=\bm{I}(\bm{X}_{i}\bm\beta \leq z\bm{1}), $$
where $\bm{X}_{ij}$, $i=1,\ldots,n$, $j=1,\ldots,n_i$ is the $j$-th row of the fixed effects design matrix  and $\bm{1}=(1,\ldots,1)^T$ is a $n_i$-vector of ones.
\end{theorem}

Next we show that the processes $W_N^{F,m}(t)$, $W_N^{\hat{F},m}(t)$ and $W_N^{\tilde{F},m} (t)$ provide valid null approximations of the process $W_N^F(t)$, see \citet{Pan05} for an alternative proof for $W_N^{\tilde{F},m} (t)$ for the special case when $M_i$ are independent standard normal random variables.

\begin{theorem}
\label{theorem1.1}
Under $H_0$ and assumptions  stated in the supplementary information,  with probability one, conditionally on data, the processes $W_N^{F,m}(t)$, $W_N^{\hat{F},m}(t)$ and $W_N^{\tilde{F},m} (t)$ converge in distribution to a zero-mean Gaussian process $G_\infty$ in the Skorokhod-space $D[-\infty,\infty]$,  where the covariance function of $G_\infty$ is
  $$K(t,s):=cov( G_{\infty}(t),G_{\infty}(s) )=\lim_{n\rightarrow \infty}K_N(t,s).$$
and where $K_N(t,s)$ is as defined in Theorem \ref{theorem1}.
\end{theorem}

 The result of Theorem \ref{theorem1.1} holds also when for $W_N^{F,m} (t)$ the residuals are ordered by fitted values in step $m$. In contrast, the result of Theorem \ref{theorem1.1} does not hold when $\bm{\Pi}_i$ is the permutation matrix since the limiting covariance function is not the same as defined in theorem \ref{theorem1}. For the same reason the result also does not hold defining
 $$\bm{y}_i^m=\bm{X}_i \hat{\bm\beta}+  \bm{\Pi}_i \bm{e}_i^P, $$
 or
 $$\bm{y}_i^m=\bm{X}_i \hat{\bm\beta}+ \bm{Z}_i\hat{\bm{b}}_i  +  \left(\bm{L}_i^T\right)^{-1}\bm{\Pi}_i\bm{L}_i^T \bm{e}_i^I, $$
 when constructing $W_N^{F,m} (t)$, nor replacing $\bm{L}_i \bm{\Pi}_i \bm{L}_i^{-1}\bm{e}_i^P$ in the definition of $W_N^{\hat{F},m} (t)$ by $\bm{\Pi}_i \bm{e}_i^P$ (see supplementary information).

We continue showing the convergence of $\tilde{W}_N^O(t)$, where  we need to additionally assume that $\bm\xi_i$, $i=1,\ldots,n$  are independent multivariate normal.  This assumption is required in order to establish that for each $i=1,\ldots,n$ and every $t$,  the transformed $\xi_i$s and $\bm\chi_i(t)$, are independent.

\begin{theorem}
\label{theormO}
Let $\bm\xi_i$, $i=1,\ldots,n$ be independent multivariate normal with $E(\bm\xi_i)=\bm{0}$ and $var(\bm\xi_i)=\bm{V}_i$, $i=1,\ldots,n$. Under $H_0$ and assumptions stated in the supplementary information,   with probability one,   the process $\tilde{W}_N^O(t )$ converges in distribution to the zero-mean Gaussian process $G_\infty$ in the Skorokhod-space $D[-\infty,\infty]$,  where the covariance function of $G_\infty$ is
  $$K(t,s):=cov( G_{\infty}(t),G_{\infty}(s) )=\lim_{n\rightarrow \infty}K_N(t,s),$$
 where
\begin{small}
\begin{eqnarray*}
K_N(t,s)&:=&\frac{1}{n}\sum_{i=1}^n   E\left(  \bm\chi_i(t,\bm{\beta},\bm{D}\bm{Z}_i^T\bm{V}_i^{-1}\bm{\xi}_i)^T  \bm{S}_i \tilde{\bm{J}}_i  \bm{V}_i \tilde{\bm{J}}_i^T  \bm{S}_i^T \bm\chi_i(s,\bm{\beta},\bm{D}\bm{Z}_i^T\bm{V}_i^{-1}\bm{\xi}_i)    \right) \\
 &-&\frac{1}{n} E\left( \left( \sum_{i=1}^n   \bm\chi_i(t,\bm{\beta},\bm{D}\bm{Z}_i^T\bm{V}_i^{-1}\bm{\xi}_i)^T\bm{S}_i\tilde{\bm{J}}_i\bm{X}_i \right)  \bm{H}^{-1}\left(  \sum_{i=1}^n  \bm{X}_i^T \tilde{\bm{J}}_i^T\bm{S}_i^T  \bm\chi_i(s,\bm{\beta},\bm{D}\bm{Z}_i^T\bm{V}_i^{-1}\bm{\xi}_i) \right)   \right),
\end{eqnarray*}
\end{small}
where $\tilde{\bm{J}}_i$ is the $i$th diagonal block of $\bm{J}$ defined as in (\ref{eq:Jot}) replacing $\bm{A}$ and $\bm{B}$ with block diagonal matrices with $i$th diagonal blocks equal to $\tilde{\bm{A}}_i$ and $\tilde{\bm{B}}_i$ defined in Lemma \ref{Rok}.
 \end{theorem}

The next theorem shows that $\tilde{W}_N^{O,m}(t )$, $\tilde{W}_N^{\hat{O},m}(t )$ and $\tilde{W}_N^{\tilde{O},m}(t )$ provide valid null approximations of the process $\tilde{W}_N^O(t)$.

\begin{theorem}
\label{theormO1}
Let $\bm\xi_i$, $i=1,\ldots,n$ be independent multivariate normal with $E(\bm\xi_i)=\bm{0}$ and $var(\bm\xi_i)=\bm{V}_i$, $i=1,\ldots,n$. Under $H_0$ and assumptions stated in the supplementary information,  with probability one,  conditionally on the data,  the processes $\tilde{W}_N^{O,m}(t )$, $\tilde{W}_N^{\hat{O},m}(t )$ and $\tilde{W}_N^{\tilde{O},m}(t )$ converge in distribution to the zero-mean Gaussian process $G_\infty$ in the Skorokhod-space $D[-\infty,\infty]$, where the covariance function of $G_\infty$ is the same as in Theorem \ref{theormO}.  
\end{theorem}



\section{Consistency under some alternative hypotheses}
\label{chap:cons}

Here we show that the tests based on $\tilde{W}_N^O(t)$  and $W_N^{F}(t)$ are consistent under some alternative hypotheses. Consistency of $W_N^O(t)$ can be established from the consistency of $\tilde{W}_N^O(t)$. Consistency of $W_N^{F^s}(t)$ can be established similarly as for $W_N^{F}(t)$ by modifying the definition of $\bm{\chi}_i(t,\bm\beta,\bm{b}_i)$.  Here we assume that $\hat{\bm{D}}$ and $\hat{\sigma}^2$ are some consistent estimators and that $n_i$, $i=1,\ldots,n$ are fixed and finite such that the model is identifiable for $\bm{D}$ and $\sigma^2$. Some assumptions are again made in order to establish that the processes have, under particular alternative hypotheses, in the limit, as $n\rightarrow\infty$, continuous sample paths (see the proofs given in the supplementary information).


We consider three alternative hypotheses which are formally stated and discussed in the supplementary information. Informally, the first alternative hypothesis (A1) states that the random effects design matrix is correctly specified while the fixed effects design matrix is not. The second alternative hypothesis (A2) states that the random effects design matrix is mis-specified whereas the fixed effects design matrix is not. Finally, the third alternative hypothesis (A3) states that both design matrices are mis-specified.
First consider the process $W_N^F(t)$. 

\begin{proposition}
\label{Alt.prop.1}
Under alternative hypotheses (A1) and (A3) there exists some $t$ such that the process $\frac{1}{\sqrt{n}}W_N^F(t)$ converges in probability towards some non-zero constant $c\neq 0$. Under alternative hypothesis (A2), the process $\frac{1}{\sqrt{n}}W_N^F(t)$ converges in probability towards zero for each $t$.
\end{proposition}

According to  Proposition \ref{Alt.prop.1}, the tests based on $W_N^F(t)$ will be powerful against alternative hypotheses (A1) and (A3) but not against alternative hypothesis (A2). Put differently, the tests based on $W_N^F(t)$ are sensitive to mis-specification of the fixed effects design matrix and robust against mis-specification of the random effects design matrix.

Assuming multivariate normality of $\bm{\xi}_i$, $i=1,\ldots,n$, the following results can be proved for the process $ \tilde{W}_N^{O}(t)$.

\begin{proposition}
\label{Alt.prop.2}
There exists some $t$ such that, under (A1) and (A2), the process $\frac{1}{\sqrt{n}}\tilde{W}_N^O(t)$ converges in probability towards some non-zero constant $c$.
\end{proposition}

According to Proposition \ref{Alt.prop.2} the tests based on $ \tilde{W}_N^{O}(t)$ will be powerful against (A1) and (A2). We show in the supplementary material that, under (A3), the result of Proposition \ref{Alt.prop.2} cannot be established for the general case. However, this will pose no practical issues, since under (A3), the model mis-specification will be, based on Proposition \ref{Alt.prop.1}, detected by inspecting $W_N^F(t)$. After correcting for the lack-of-fit due to mis-specifying the fixed effects design matrix, this will then be detected inspecting $ \tilde{W}_N^{O}(t)$ (case (ii) of Proposition \ref{Alt.prop.2}).



\section{Simulation results}
\label{sec:sim}
Here we show a set of selected simulation results. Complete simulation results, including also the same simulation design as in \citet{Pan05}, are shown in supplementary information. The outcome was simulated from
$$ y_{ij}= - 1 +0.25 X_{ij,1}+0.5X_{ij,2}+\beta_3X_{ij,1}^2+b_{i,0}+b_{i,1}X_{ij,1}+\epsilon_{ij}\mbox{, }j=1,\ldots,n_i\mbox{, }i=1,\ldots,n,$$
simulating the same number of observations in each cluster, where $X_{ij,1}\sim U(0,1)$, $X_{ij,2}\sim U(0,1)$, $\epsilon_{ij}\sim N(0,0.5)$, $b_{i,0}\sim N(0,0.25)$ and $b_{i,1}\sim N(0,\sigma_{b,1}^2)$. $b_{i,0}$ and $b_{i,1}$ were simulated independently and $\bm{\epsilon}_i=(\epsilon_{i1},\ldots,\epsilon_{in_i})^T$ and $\bm{b}_i=(b_{i,0},b_{i,1})^T$ were also independent.

The analysis was performed in R (R version 3.4.3,\citet{R}) using the R package \textbf{gofLMM} (available on GitHub, \url{rokblagus/gofLMM}). The LMMs were fitted by using the function \textbf{lme} from the \textbf{nlme} package \citep{BatesBook}. The variance parameters were estimated by REML. For each simulated scenario we calculated the KS and CvM type test statistics based on processes $\tilde{W}_N^O(t)$  and $W_N^F(t)$ 
and the $p$-values were estimated by using $M=500$ random simulations/sign-flips  as described in section \ref{gof}.  Each step of the simulation was repeated $5000$ times; the simulation margin of errors are $\pm 0.003$, $\pm 0.006$ and $\pm 0.008$ for $\alpha=0.01$, 0.05 and $0.1$, respectively. 

Throughout, the results when using \textit{individual} or \textit{cluster} residuals were very similar, hence only the results for \textit{individual}   residuals are shown here. We only show the results for the CvM type test statistic since, in general, using CvM type test statistic was more powerful (see supplementary information).

\subsection{Example I - size under normal errors and random effects}
\label{ExampleI}

Here the outcome is simulated using $\sigma_{b,1}^2=0.25$, $\beta_3=0$, $n=50,75$ and for each $n$, $n_i=5,10,20$. The fitted model was the same as the simulated model.


With large $n$ and/or $n_i$, the empirical size of the tests were close to nominal levels (Table \ref{Example1tab}). With small $n$, the simulation approach was conservative, while our novel simulation approach using Cholesky decomposition was too liberal (more obvious when testing the goodness-of-fit for the entire LMM). Same behavior was observed, although too a much lesser extent, for the proposed approach based on sign-flipping. We find this a consequence of using an asymptotic approximation when calculating $\bm{e}_i^C$, since this was not observed when using matrices $\bm{A}$ and $\bm{B}$ when calculating $\bm{e}^C$ (see supplementary information). Overall, the performance of the tests, especially the one using the proposed sign-flipping approach, was satisfactory.

\begin{table}[ht]
\caption{\small{Empirical sizes of the cusum tests (O - test for the entire model, F - test for the fixed effects part of the model) using  CvM type test statistics using different procedures to approximate the null distribution (Pan - the simulation approach, Sim - the proposed simulation approach using Cholesky decomposition, SF - the proposed approach based on sign-flipping).}}
\label{Example1tab}
\spacingset{1}
\centering
\begingroup\scriptsize
\begin{tabular}{rrrrrrrrrr}
  \hline
Example&$\alpha$ & $n$ & $n_i$ & Pan:O & Pan:F & Sim:O & Sim:F & SF:O & SF:F \\
  \hline
I&0.10 & 50 & 5 & 0.0816 & 0.0934 & 0.1268 & 0.0980 & 0.1172 & 0.0978 \\
&   &  & 10 & 0.0862 & 0.0862 & 0.1138 & 0.0908 & 0.1072 & 0.0922 \\
&   &  & 20 & 0.0836 & 0.0900 & 0.1054 & 0.1000 & 0.1038 & 0.0998 \\
 &  & 75 & 5 & 0.0878 & 0.0996 & 0.1332 & 0.1046 & 0.1206 & 0.1032 \\
  & &  & 10 & 0.0986 & 0.0982 & 0.1258 & 0.0996 & 0.1202 & 0.1002 \\
&   &  & 20 & 0.0882 & 0.0934 & 0.1068 & 0.0948 & 0.1068 & 0.0952 \\\cline{2-10}
&  0.05 & 50 & 5 & 0.0306 & 0.0396 & 0.0672 & 0.0466 & 0.0546 & 0.0454 \\
&   &  & 10 & 0.0354 & 0.0372 & 0.0618 & 0.0430 & 0.0542 & 0.0428 \\
 &  &  & 20 & 0.0344 & 0.0430 & 0.0536 & 0.0494 & 0.0530 & 0.0472 \\
 &  & 75 & 5 & 0.0370 & 0.0518 & 0.0688 & 0.0522 & 0.0616 & 0.0532 \\
 &  &  & 10 & 0.0424 & 0.0454 & 0.0634 & 0.0498 & 0.0598 & 0.0502 \\
 &  &  & 20 & 0.0400 & 0.0416 & 0.0528 & 0.0448 & 0.0538 & 0.0450 \\\cline{2-10}
 & 0.01 & 50 & 5 & 0.0042 & 0.0052 & 0.0112 & 0.0074 & 0.0124 & 0.0080 \\
 &  &  & 10 & 0.0034 & 0.0056 & 0.0110 & 0.0070 & 0.0104 & 0.0080 \\
 &  &  & 20 & 0.0050 & 0.0074 & 0.0110 & 0.0098 & 0.0130 & 0.0096 \\
 &  & 75 & 5 & 0.0052 & 0.0100 & 0.0134 & 0.0104 & 0.0112 & 0.0124 \\
 &  &  & 10 & 0.0056 & 0.0102 & 0.0134 & 0.0108 & 0.0144 & 0.0106 \\
 &  &  & 20 & 0.0082 & 0.0074 & 0.0114 & 0.0090 & 0.0098 & 0.0094 \\
   \hline
 II&0.10 & 50 & 5 & 0.0418 & 0.0872 & 0.1036 & 0.0962 & 0.0954 & 0.0944 \\
  & &  & 10 & 0.0452 & 0.0892 & 0.0984 & 0.0944 & 0.1020 & 0.0954 \\
  & &  & 20 & 0.0501 & 0.0931 & 0.0879 & 0.0959 & 0.0953 & 0.0993 \\
  & & 75 & 5 & 0.0485 & 0.0958 & 0.1218 & 0.1036 & 0.1066 & 0.1020 \\
  & &  & 10 & 0.0558 & 0.0864 & 0.1097 & 0.0906 & 0.1174 & 0.0916 \\
  & &  & 20 & 0.0670 & 0.0870 & 0.1012 & 0.0864 & 0.1118 & 0.0880 \\\cline{2-10}
  &0.05 & 50 & 5 & 0.0102 & 0.0346 & 0.0440 & 0.0424 & 0.0446 & 0.0468 \\
  & &  & 10 & 0.0152 & 0.0408 & 0.0424 & 0.0492 & 0.0464 & 0.0486 \\
  & &  & 20 & 0.0158 & 0.0430 & 0.0400 & 0.0472 & 0.0488 & 0.0486 \\
  & & 75 & 5 & 0.0148 & 0.0457 & 0.0525 & 0.0471 & 0.0547 & 0.0501 \\
  & &  & 10 & 0.0203 & 0.0385 & 0.0493 & 0.0427 & 0.0578 & 0.0445 \\
  & &  & 20 & 0.0224 & 0.0426 & 0.0478 & 0.0476 & 0.0548 & 0.0470 \\\cline{2-10}
  &0.01 & 50 & 5 & 0.0000 & 0.0052 & 0.0046 & 0.0086 & 0.0080 & 0.0096 \\
  & &  & 10 & 0.0008 & 0.0074 & 0.0058 & 0.0096 & 0.0076 & 0.0092 \\
  & &  & 20 & 0.0014 & 0.0062 & 0.0066 & 0.0098 & 0.0092 & 0.0102 \\
  & & 75 & 5 & 0.0010 & 0.0052 & 0.0060 & 0.0070 & 0.0084 & 0.0088 \\
  & &  & 10 & 0.0014 & 0.0056 & 0.0087 & 0.0083 & 0.0119 & 0.0097 \\
  & &  & 20 & 0.0016 & 0.0058 & 0.0078 & 0.0098 & 0.0106 & 0.0082 \\

 \hline
\end{tabular}
\endgroup
\end{table}

\spacingset{1}

\subsection{Example II - non-normal errors and random effects}
\label{ExampleIa}

Here the outcome was simulated as in section \ref{ExampleI}, but simulating $\epsilon_{ij}$, $b_{i,0}$ and $b_{i,1}$ independently from a zero mean gamma distribution with parameters shape and scale set to 1 and 2, respectively. The fitted model was the same as the simulated model.

\begin{figure}[h!]
\caption{Empirical sizes of the cusum tests (O - test for the entire model, F - test for the fixed effects part of the model) using  CvM type test statistics using different procedures to approximate the null distribution (simulation.Pan - the simulation approach, simulation - the proposed simulation approach using Cholesky decomposition, sign-flip - the proposed approach based on sign-flipping) for the example with non-normal errors and random effects.}
\label{Example1a}
\centering
 \resizebox{160mm}{!}{\includegraphics{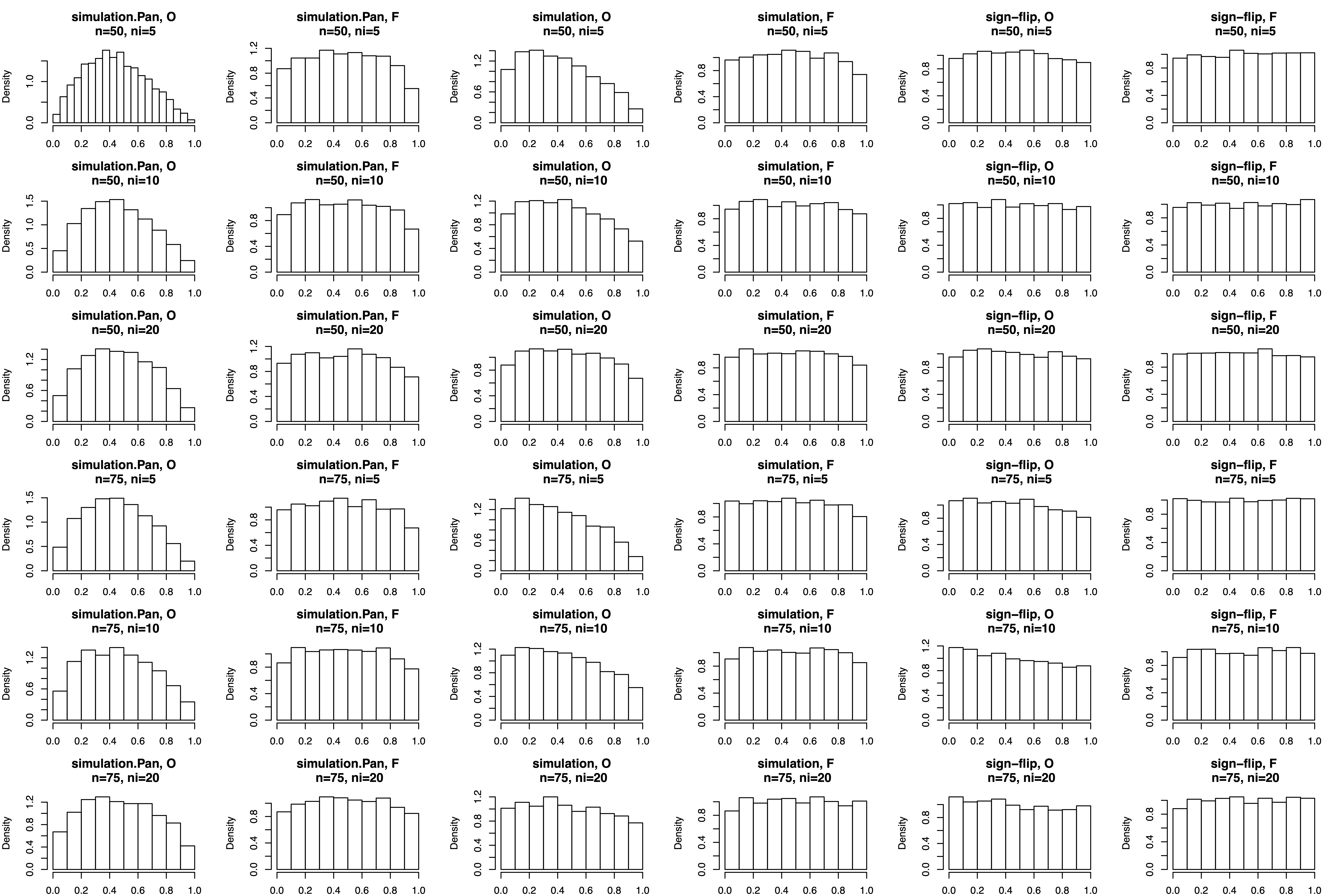}}
\end{figure}

In this example the tests based on the simulation approach did not perform well obtaining distributions of the estimated $p$-values which were not uniform (see Figure \ref{Example1a} and Table \ref{Example1tab}). In contrast, the proposed approach based on sign-flipping performed very similarly as in  Example I, showing its robustness against non-normal random effects and errors.

 \clearpage
\newpage

\subsection{Example III - mis-specified random effects design matrix }
\label{ExampleII}

The outcome was simulated as presented in section \ref{sec:sim}, using $\sigma_{b,1}^2=0.5,1,1.5$, $\beta_3=0$, $n=50,75$ and for each $n$, $n_i=10$. The fixed effects part of the fitted model was correctly specified, but the random effects part included only random intercept.


As suggested by our theoretical results, the empirical sizes of the tests for the fixed effects parts of the model were similar to nominal levels, demonstrating their robustness against mis-specification of the random effects design matrix. The tests for the overall goodness-of-fit of the model rejected the null hypothesis more often then the nominal level. Rejection rates were larger with larger $n$ and/or $\sigma_{b,1}^2$, with the proposed approach based on sign-flipping being the most powerful (Table \ref{Example2tab}).

\begin{table}[h!]
\caption{Empirical powers/sizes of the cusum tests (O - test for the entire model, F - test for the fixed effects part of the model) using  CvM type test statistics using different procedures to approximate the null distribution ( Pan - the simulation approach, Sim - the proposed simulation approach using Cholesky decomposition, SF - the proposed approach based on sign-flipping) for the examples with mis-specified random effects and fixed effects design matrices.}
\label{Example2tab}
\spacingset{1}
\centering
\begingroup\scriptsize
\begin{tabular}{rrrrrrrrrr}
  \hline
Example&$\alpha$ & $n$ & $\sigma_{b,1}^2$/$\beta_3$ & Pan:O & Pan:F & Sim:O & Sim:F & SF:O & SF:F \\
  \hline
III&0.10 & 50 & 0.5 & 0.1520 & 0.0948 & 0.2194 & 0.1010 & 0.2210 & 0.1010 \\
  & &  & 1.0 & 0.2314 & 0.0946 & 0.3604 & 0.1006 & 0.3642 & 0.1006 \\
  & &  & 1.5 & 0.3374 & 0.0904 & 0.4990 & 0.0986 & 0.5050 & 0.0970 \\
  & & 75 & 0.5 & 0.1812 & 0.0938 & 0.2544 & 0.0974 & 0.2582 & 0.0976 \\
  & &  & 1.0 & 0.3204 & 0.0990 & 0.4368 & 0.1048 & 0.4418 & 0.1058 \\
  & &  & 1.5 & 0.4418 & 0.1036 & 0.5766 & 0.1034 & 0.5832 & 0.1034 \\ \cline{2-10}
  &0.05 & 50 & 0.5 & 0.0734 & 0.0470 & 0.1378 & 0.0522 & 0.1428 & 0.0532 \\
  & &  & 1.0 & 0.1302 & 0.0444 & 0.2662 & 0.0498 & 0.2720 & 0.0490 \\
  & &  & 1.5 & 0.2168 & 0.0418 & 0.3892 & 0.0476 & 0.4028 & 0.0478 \\
  & & 75 & 0.5 & 0.0974 & 0.0444 & 0.1630 & 0.0428 & 0.1650 & 0.0438 \\
  & &  & 1.0 & 0.2050 & 0.0478 & 0.3424 & 0.0528 & 0.3406 & 0.0522 \\
  & &  & 1.5 & 0.3200 & 0.0492 & 0.4786 & 0.0502 & 0.4840 & 0.0528 \\ \cline{2-10}
  &0.01 & 50 & 0.5 & 0.0132 & 0.0074 & 0.0476 & 0.0092 & 0.0486 & 0.0086 \\
  & &  & 1.0 & 0.0340 & 0.0064 & 0.1196 & 0.0094 & 0.1302 & 0.0098 \\
  & &  & 1.5 & 0.0680 & 0.0056 & 0.2258 & 0.0072 & 0.2438 & 0.0074 \\
  & & 75 & 0.5 & 0.0194 & 0.0066 & 0.0510 & 0.0094 & 0.0568 & 0.0084 \\
  & &  & 1.0 & 0.0616 & 0.0090 & 0.1642 & 0.0098 & 0.1756 & 0.0094 \\
  & &  & 1.5 & 0.1350 & 0.0078 & 0.3122 & 0.0108 & 0.3216 & 0.0100 \\
     \hline
  IV&   0.10 & 50 & 0.5 & 0.0866 & 0.1358 & 0.1218 & 0.1464 & 0.1060 & 0.1470 \\
  & &  & 1.0 & 0.1102 & 0.4162 & 0.1482 & 0.4392 & 0.1260 & 0.4372 \\
  & &  & 1.5 & 0.2600 & 0.8038 & 0.3290 & 0.8180 & 0.3050 & 0.8158 \\
  & & 75 & 0.5 & 0.0964 & 0.1504 & 0.1204 & 0.1550 & 0.1106 & 0.1560 \\
  & &  & 1.0 & 0.1458 & 0.5664 & 0.1840 & 0.5778 & 0.1676 & 0.5734 \\
  & &  & 1.5 & 0.3766 & 0.9274 & 0.4460 & 0.9360 & 0.4248 & 0.9346 \\ \cline{2-10}
  &0.05 & 50 & 0.5 & 0.0380 & 0.0736 & 0.0608 & 0.0810 & 0.0480 & 0.0810 \\
  & &  & 1.0 & 0.0530 & 0.2968 & 0.0788 & 0.3242 & 0.0712 & 0.3226 \\
  & &  & 1.5 & 0.1482 & 0.6980 & 0.2186 & 0.7288 & 0.1900 & 0.7312 \\
  & & 75 & 0.5 & 0.0472 & 0.0858 & 0.0624 & 0.0890 & 0.0592 & 0.0906 \\
  & &  & 1.0 & 0.0682 & 0.4482 & 0.1060 & 0.4660 & 0.0954 & 0.4658 \\
  & &  & 1.5 & 0.2466 & 0.8762 & 0.3216 & 0.8902 & 0.3028 & 0.8910 \\ \cline{2-10}
  &0.01 & 50 & 0.5 & 0.0044 & 0.0150 & 0.0140 & 0.0204 & 0.0092 & 0.0200 \\
  & &  & 1.0 & 0.0070 & 0.1112 & 0.0190 & 0.1444 & 0.0154 & 0.1436 \\
  & &  & 1.5 & 0.0332 & 0.4202 & 0.0774 & 0.4944 & 0.0608 & 0.5064 \\
  & & 75 & 0.5 & 0.0076 & 0.0196 & 0.0154 & 0.0242 & 0.0126 & 0.0232 \\
  & &  & 1.0 & 0.0130 & 0.2080 & 0.0264 & 0.2428 & 0.0212 & 0.2458 \\
  & &  & 1.5 & 0.0756 & 0.6952 & 0.1310 & 0.7376 & 0.1156 & 0.7438 \\
   \hline

\end{tabular}
\endgroup
\end{table}

\spacingset{1}

\subsection{Example IV - mis-specified fixed effects design matrix}
\label{ExampleIII}

The outcome was simulated as presented in section \ref{sec:sim}, using $\sigma_{b,1}^2=0.25$, $\beta_3=0.5,1,1.5$, $n=50,75$ and for each $n$, $n_i=10$. The random effects part of the fitted model was correctly specified, but the fixed effects part included only linear effects of the covariates.


Empirical rejection rates of all tests were larger than the nominal level, showing that the tests are powerful against this alternative. Not surprisingly, the rejection rates when using $W_N^F(t)$ were larger then when using $W_N^O(t)$ (Table \ref{Example2tab}).

\section{Application}
\label{real}

We apply the proposed methodology to the CD4 count data available in the R package \textbf{JSM}. The data contain longitudinal measurements for 467 patients (in total 1405 measurements) which were randomly assigned to either zalcitabine or didanosine antiretroviral treatment.

\begin{figure}[h!]
\begin{center}
\resizebox{120mm}{!}{\includegraphics{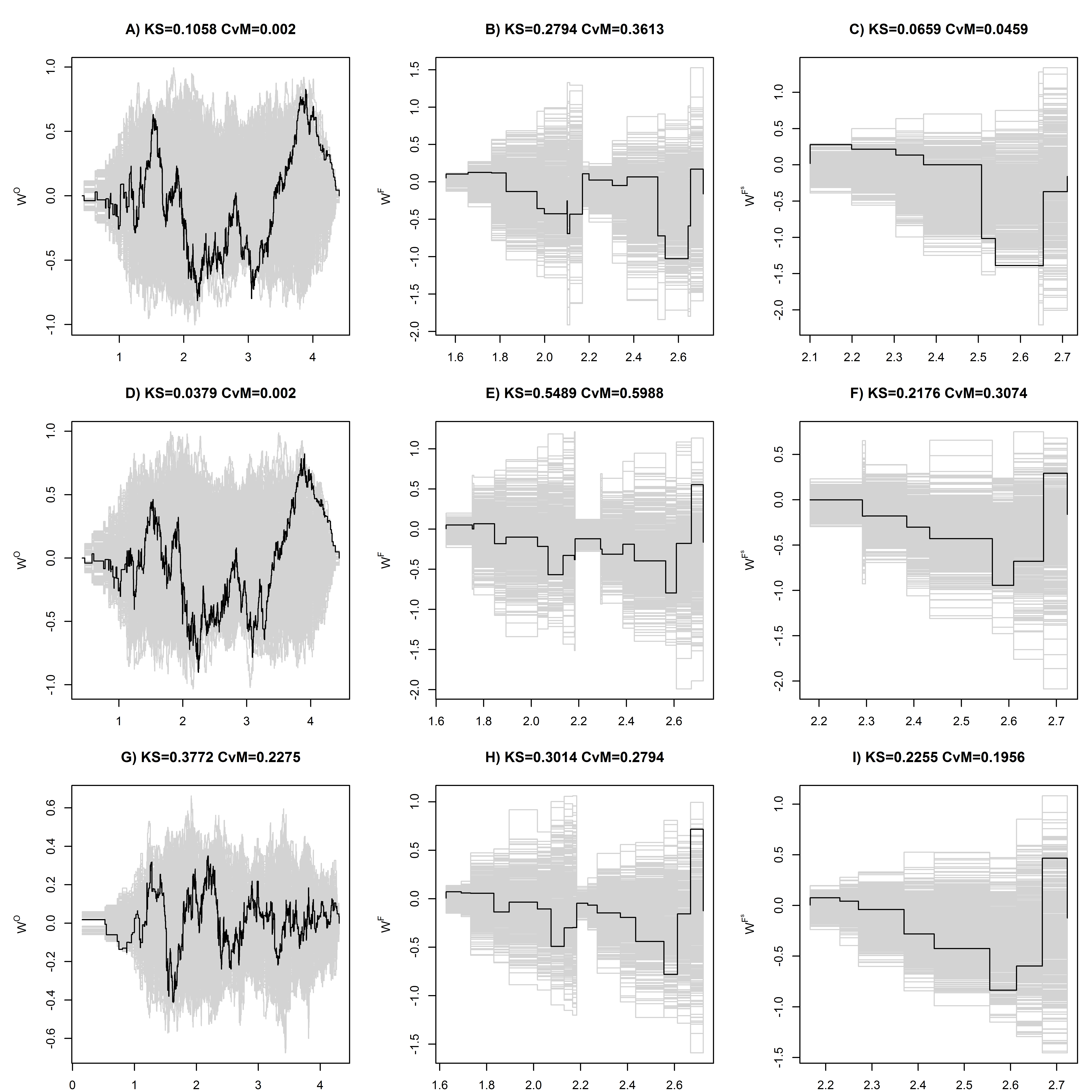}}
\caption{Cusum processes $\tilde{W}_N^O(t)$, $W_N^F(t)$ and $W_N^{F^s}(t)$ (left, center and right panels, respectively) for different models (rows 1 to 3 correspond to Models 1 to 3, respectively), for the CD4 dataset. Sign-flipping approach.}
\label{fig.real.ill}
\end{center}
\end{figure}

Let $CD4_{ij}$ denote the CD4 cell counts for the $j$th measurement of individual $i$. Let $R_{i}$ denote the treatment received by individual $i$ (1 is used for zalcitabine treatment) and let $T_{ij}$ denote the time at which the $j$th measurement for individual $i$ was taken. With $A_{i}$ we denote the AIDS indicator at the start of the study (1 is used for no AIDS). We fit the following models to the data:
\begin{equation}
\label{real.model.1}
\mbox{Model 1 }CD4_{ij}=\beta_1+\beta_2 A_{ij}+\beta_3T_{ij}+\beta_4T_{ij}R_{ij}+b_{i},
\end{equation}
\begin{equation}
\label{real.model.2}
\mbox{Model 2 }CD4_{ij}=\beta_1+\beta_2 A_{ij}+\beta_3T_{ij}+\beta_4T_{ij}R_{ij}+\beta_5T_{ij}^2+\beta_6T_{ij}^2R_{ij}+b_i,
\end{equation}
\begin{equation}
\label{real.model.3}
\mbox{Model 3 }CD4_{ij}=\beta_1+\beta_2 A_{ij}+\beta_3T_{ij}+\beta_4T_{ij}R_{ij}+\beta_5T_{ij}^2+\beta_6T_{ij}^2R_{ij}+b_{i1}+b_{i2}T_{ij}.
\end{equation}
For each fitted model we plot, with black lines, $\tilde{W}_N^O(t)$, $W_N^F(t)$ and $W_N^{F^s}(t)$, where the subset in the last process includes all the estimated coefficients associated with variable $T_{ij}$ (also those where $T_{ij}$ is included in the interaction term)  along with $\tilde{W}_N^{O,m}(t)$, $W_N^{F,m}(t)$ and $W_N^{F^s,m}(t)$  shown as gray lines where we used sign-flipping  with $M=500$. 
The results are shown in Figure \ref{fig.real.ill}.

Inspection of the process $W_N^{F^s}(t)$ for Model 1 reveals that a quadratic term is missing for $T_{ij}$ (Figure \ref{fig.real.ill}, row 1, column 3), and we can see that including the quadratic term (Model 2) has substantially improved the fit of the model (Figure \ref{fig.real.ill}, row 2, columns 2 and 3). The $p$-value obtained from the process $\tilde{W}_N^O(t)$ is still significant at the 0.05 level (Figure \ref{fig.real.ill}, row 2, column 1), suggesting that the random effects structure is misspecified. Including also the random slope for $T_{ij}$ substantially improves the model's fit (Figure \ref{fig.real.ill}, row 3, column 1).

\section{Discussion and Conclusions}
\label{disc}

We showed how cusum processes can be used to test the assumed functional form of fitted LMMs.  We proposed a procedure based on inspecting plots of the cusum processes ordering the transformed residuals by the \textit{individual} predicted values. Ordering the residual by \textit{cluster} predicted values on the other hand allows investigating the assumed functional form only of the assumed fixed effects design matrix (or some subset thereof). We showed that when appropriately transforming the residuals, the cusum process for the entire model is, with a reasonably large $n$, expected to fluctuate around zero when the fixed and random effects design matrices are correctly specified, while it is not when either (or both) fixed and random effects design matrices are mis-specified. In contrast, the cusum process targeting the fixed effects part of the LMM was shown to fluctuate around zero when the fixed effects design matrix is correctly specified, regardless of the (in)correct specification of the random effects design matrix.

Observed fluctuations/deviations from zero can be compared and evaluated by means of $p$-values by using sign-flipping/bootstrap or (novel) simulation approach. The later is appealing with a large sample size and/or complex models since, in contrast to the sign-flipping/bootstrap approach, does not require re-estimating the LMM. It was shown theoretically that all approaches are asymptotically equivalent, however it was demonstrated by means of a large Monte-Carlo simulation study, that the approach based on sign-flipping yields better results with smaller sample size and/or non-normally distributed random effects and/or errors in terms of size while obtaining similar power. The proposed simulation approach was more powerful than the simulation approach proposed by \citet{Pan05}, however both approaches require a sufficiently large sample size, especially with non-normal random effects and/or errors,  to achieve valid inference under $H_0$.

To prove the consistency of the test for the entire LMM we assumed that random effects and errors are multivariate normal, while this was not required when inspecting only the fixed effects design matrix. In the simulation study we observed that, even with a small sample size, the type I errors were not inflated with non-normal random errors and effects when using the proposed sign-flipping/bootstrap approach. 

Our theoretical results imply that the processes can be based on either  \textit{individual} or \textit{cluster} residuals. 
In the simulation study we observed that using either option leads to very similar results. In our reanalysis of the CD4 count data from the R's \textbf{JSM} package we noticed, however that using \textit{cluster} or \textit{individual} residuals can lead to contrasting conclusions (supplementary information). While the finite sample performance of the proposed methodology with unbalanced design seems an interesting area for further research, we did not address this in more detail.

The proposed tests are non-directional, hence potentially being less powerful then directional tests. However, the approach proposed here can easily be combined/supplemented by using directional tests, where the plots of the cusum processes can be used as a guide towards (the potentially more powerful) directional tests. E.g. in our simulated example, Wald test for the quadratic term of the fixed effect covariate could also be used to test if a quadratic association (implied by inspecting the plot of the cusum process) needs to be modeled, although realizing that this particular functional form for this particular covariate is suitable, would be more challenging without investigating the plot of the cusum process.

The proposed methods were implemented in R package \textbf{gofLMM} utilizing the \textbf{nlme} package \citep{BatesR}, adopting its great flexibility. At the time of preparing this manuscript the R package was still at the development stage, but we are planing to make the package available on CRAN in the near future.

While we in detail considered only single-level LMMs, we showed that the proposed methodology could easily be adapted to multiple, nested levels of random effects, but at a cost of notational inconvenience. In principle the methodology presented here could be extended to GLMMs. However, further extensions to nonlinear link GLMMs could be problematic when trying to distinguish between the reasons for the (possible) lack-of-fit due to fixed or random effect design matrices.

\begin{small}
\subsection*{Supplementary information}
The  on-line  supplementary  information  document  contains  proofs  of  the  theorems  and  other technical details. Complete simulation results are also shown.

\begin{itemize}
\item \href{https://www.dropbox.com/s/ehedeizhblfa26w/supplementaryMaterialV4.pdf?dl=0}{Supplementary
material (link)}
\end{itemize}

\subsection*{Acknowledgements}
JP is a young researcher funded by the Slovenian Research Agency (ARRS). RB acknowledges the financial support by ARRS (Predicting rare events more accurately, N1-0035; Methodology for data analysis in medical sciences, P3-0154).

\end{small}

\spacingset{0.1} 
\begin{footnotesize}
\bibliographystyle{plainnat}
\bibliography{goflmm}
\end{footnotesize}


\end{document}